\documentclass[12pt,reqno]{amsart}
\usepackage{amscd,amssymb,amsmath}
\usepackage[arrow,matrix,curve]{xy}
\newtheorem{theorem}{Theorem}[section]
\newtheorem{proposition}[theorem]{Proposition}
\newtheorem{lemma}[theorem]{Lemma}
\newtheorem{corollary}[theorem]{Corollary}

\theoremstyle{definition}
\newtheorem{definition}[theorem]{Definition}

\newtheorem{remark}[theorem]{Remark}

\newtheorem{question}[theorem]{Question}
\newtheorem{hypothesis}[theorem]{Hypothesis}

\numberwithin{equation}{section}
\numberwithin{figure}{section}

\def\Z {{\mathbb Z}}

\def\N {{\mathbb N}}

\def\R {{\mathbb R}}
\def\C {{\mathbb C}}

\def\H {{\mathcal{H}}} \def\P {{\mathcal{P}}} \def\M {{\mathcal{M}}}
 \def\G {{\mathcal{G}}}

\def\ie {{\it i.e.\ }}

\def\cf {\hbox{\it cf.\ }}

\def\CP {{\C P}}

\def\vol{{\it vol}}
 
\def\area{{\it area}} 
\def\sys{{\it sys}} 
 
\def\stsys{{\it stsys}} 
\def\confsys{{\it confsys}}

\def\CT{{\it CT}}
\def\SR{{\it SR}}
\def\PD{{\it PD}}

\def\gen{{s}}
\def\v{{v}}

\normalbaselines 

\input epsf 		%%for the figure

\begin{document}
%\large

\author[M.~Katz]{Mikhail Katz$^{*}$} \address{Department of
Mathematics and Statistics, Bar Ilan University, Ramat Gan 52900
Israel} \email{katzmik@math.biu.ac.il} \thanks{$^{*}$Supported by the
Israel Science Foundation (grant no.\ 620/00-10.0).  Conditionally
accepted for publication in $\underline{\hbox{Commentarii Mathematici
Helvetici}}$.}

\title[
Four-manifold systoles and surjectivity of period map
]%
{
Four-manifold systoles and surjectivity of period map
}

\begin{abstract}
P. Buser and P. Sarnak showed in 1994 that the maximum, over the
moduli space of Riemann surfaces of genus $\gen$, of the least
conformal length of a nonseparating loop, is logarithmic in $\gen$.
We present an application of (polynomially) dense Euclidean packings,
to estimates for an analogous 2-dimensional conformal systolic
invariant of a 4-manifold $X$ with indefinite intersection form.  The
estimate turns out to be polynomial, rather than logarithmic, in
$\chi(X)$, if the conjectured surjectivity of the period map is
correct.  Such surjectivity is targeted by the current work in gauge
theory.  The surjectivity allows one to insert suitable lattices with
metric properties prescribed in advance, into the second de~Rham
cohomology group of $X$, as its integer lattice.  The idea is to adapt
the well-known Lorentzian construction of the Leech lattice, by
replacing the Leech lattice by the Conway-Thompson unimodular lattices
which define asymptotically dense packings.  The final step can be
described, in terms of the successive minima $\lambda_i$, as deforming
a $\lambda_2$-bound into a $\lambda_1$-bound, illustrated by
Figure~\ref{101}.
\end{abstract}

\maketitle

\tableofcontents

\section{Schottky problem, surjectivity conjecture, and main theorem}

The work of P. Buser and P. Sarnak \cite{BS} on Riemann surfaces in
connection with the Schottky problem shows that the maximum, over the
moduli space, of the least conformal length of a nonseparating loop
behaves {\it logarithmically} in the genus, \cf \eqref{bs} and also
M.~Gromov's result \cite[Theorem $5.5.C'$]{Gro1}.

We provide a lower bound which is {\it polynomial} in the second Betti
number, for the analogous 2-dimensional conformal systolic invariant
for a 4-manifold $X$ with indefinite intersection form, modulo the
conjectured surjectivity of the period map, targeted in the current
work \cite{ADK}.  Our bound currently depends on such surjectivity,
but see \ref{(7.2)}.  In the case $b^+(X) = 1$ targeted in \cite{ADK},
such surjectivity is expressed in Hypothesis~\ref{(1.1)}.

Let $(X,g)$ be a Riemannian 4-manifold.  Let $*: H^2_{dR}(X) \to
H^2_{dR}(X) $ be the Hodge star operator in de Rham cohomology
identified with the space $\H$ of harmonic 2-forms on $X$.  Assume
that $b^+(X) = 1$, so that the selfdual subspace (\ie the
$(+1)$-eigenspace of the Hodge star operator) is 1-dimensional.
Recall that the cup-product form in $H^2(X)$ is dual to the
intersection form in $H_2(X)$.

\begin{hypothesis}
\label{(1.1)}  
For every line $V$ in the positive cone in $H^2_{dR}(X)$ defined by
the cup product form, there is a metric $g$ on $X$ whose selfdual
subspace is exactly $V$.
\end{hypothesis}

Given a lattice $L$ equipped with a norm $\|\;\|$, we denote by
\begin{equation}
\label{lambda}
\lambda_1(L)= \lambda_1(L,\|\;\|)
\end{equation}
the least norm of a nonzero lattice vector.  The $\lambda_1$ notation
fits in with the successive minima $\lambda_i$ of a lattice, studied
in lattice theory, \cf \cite[p.~58]{GruL}, \cite[Section~4]{BanK}, and
Definition~\ref{(2.three)}.

\begin{theorem}
\label{T}  
Let $n\in \N$ and consider the complex projective plane blown up at
$n$ points, $P_n=\CP^2\#n \overline{\CP}^2,$ where bar denotes
reversal of orientation, while $\#$ is connected sum.  Assume that
Hypothesis~\ref{(1.1)} is satisfied for such manifolds.  Then
\begin{equation}
\label{(3.3)}
C^{-1} \sqrt{n} < \sup_g \left\{ \lambda_1\left( H^2(P_n,\Z),
|\;|^{\phantom{A}}_{L^2} \right) \right\} ^2 < C n, \;\forall n>0,
\end{equation}
where $C>0$ is a numerical constant, the supremum is over all smooth
metrics $g$ on $\CP^2\#n \overline{\CP}^2$, and
$|\;|_{L^2}^{\phantom{a}}$ is the norm \eqref{(1.2)} defined by $g$.
\end{theorem}
Here the upper bound may be replaced by $\frac{2}{3}(n+1)$ by the
estimate~\eqref{(3.1)}, while the lower bound, by $\sqrt{k(n)}$, where
$k(n)$ is asymptotic to $\frac{n}{2\pi e}$ as $n\to \infty$, \cf
Theorem~\ref{(1.4)}.  The theorem is proved in
Section~\ref{punchline}.  The desired metric is specified in
formula~\eqref{F} in terms of inversion of the period map.

A number of systolic inequalities are now available in the literature.
Nontrivial cup product relations lead to stable systolic inequalities
\cite{BanK} (\cf inequality~\eqref{(3.1)} below), some of them sharp
\cite{BanK2,NV}.  Meanwhile, nontrivial Massey products also admit
systolic repercussions, \cf [KKS].  For ordinary (rather than stable)
systoles, systolic freedom prevails as soon as we go beyond loops.
Here ``systolic freedom'' refers to the absence of systolic
inequalities, \ie the existence of sequences of metrics violating such
potential inequalities.  Such a phenomenon for the middle dimensional
systole was first described in detail by the author in \cite{Ka1}.
M.~Gromov's original seminal example is described in \cite{Be} and
\cite[p.~268]{Gro3}, as well as \cite[section 4.2]{CK}.  Further
generalisations of systolic freedom were obtained in \cite{BabK},
\cite{Fr}, \cite{KS1,KS}, \cite{Bab}, \cite{Ka2}.  See the recent
survey \cite[Figure~4.1]{CK} for a 2-D map of systolic geometry, which
places such results in mutual relation.

\begin{question}
\label{(7.2)}  
Can one eliminate the dependence of our Theorem~\ref{T} on the
surjectivity conjecture?  Recent discussions with C. LeBrun and
P. Biran suggest that one may be able to remove the dependence on the
conjectured surjectivity of the period map, at least in the case of
the blow-ups of the projective plane, by exploiting the action of the
automorphism group of the intersection form, \cf Lemma~\ref{(8.1)} and
Remark~\ref{kotschick}.  This would work if one can show the existence
of metrics adapted to symplectic forms which represent classes from a
suitable fundamental domain for the action, \cf \cite[Theorem 3.2]{Bi}
and Remark~\ref{(3.4)}.
\end{question}

\begin{question}
\label{(4.5)}  
Can one improve the lower bound in \eqref{(3.3)} to {\em linear\/}
dependence on $n$?  Here one could envision an averaging argument,
using Siegel's formula as in \cite[Theorem~9.5]{MH}, over integral
vectors satisfying $q_{n,1}(v)=-p$.  Here one seeks a vector $v\in
\R^{n,1}$ such that the integer lattice $\Z^{n,1}\subset \R^{n,1}$ has
the Conway-Thompson behavior \eqref{71} with respect to the positive
definite form $\SR(q_{n,1},v)$.
\end{question}

\begin{question}
\label{(4.6)}  
Is there an {\em asymptotically infinite} lower bound similar to
Theorem~\ref{T} for the stable 2-systole in place of the conformal
2-systole?  This is related to understanding the discrepancy between
the comass norm and the $L^2$ norm in 2-dimensional cohomology.  Note
that Remark~\ref{(3.four)}, concerning the 1-systoles of surfaces,
suggests that {\it a priori} there may exist, instead, an {\em
asymptotically vanishing} upper bound for the stable 2-systole, \cf
\eqref{42} in the definite case.
\end{question}

The present work is organized as follows.  Section~\ref{(3)}
introduces the Conway-Thompson lattices and describes the idea of the
proof.  Section~\ref{(two)} defines the $L^2$-norm in cohomology,
describes its relation to the intersection form, and discusses the
successive minima of a lattice.  Section~\ref{(three)} defines the
conformal and stable systoles.  Section~\ref{four} discusses the
definite case.  Note that our main Theorem~\ref{T} can be thought of
as a higher-dimensional analogue of the Buser-Sarnak theorem,
presented in Section~\ref{(2)}.  Section~\ref{six} explains a useful
sign reversal relation between definite and indefinite forms.
Section~\ref{8} describes a Lorentzian construction of lattices
inspired by a result of J.~Conway and N.~Sloane, and presents a lower
bound for the second successive minimum.  Section~\ref{(5)} presents
the necessary linear algebraic ingredient.  Section~\ref{(7)} deforms
a lower bound for the second successive minimum, into a lower bound
for the first successive minimum.  The proof is completed by a
successive minimum calculation in Section~\ref{punchline}.

\section{Conway-Thompson 
lattices $\CT_n$ and idea of proof}
\label{(3)}
The surjectivity of the period map (see Hypothesis \ref{(1.1)})
furnishes a lot of latitude in prescribing the position of the integer
lattice in middle-dimensional de Rham cohomology, with respect to the
$L^2$-norm.  In particular, we show that the least norm,
$\lambda_1\left(H^p(X,\Z)_\R, |\;|^{\phantom{a}}_{L^2}\right)$, of a
nonzero lattice element, can be made arbitrarily large as the Betti
number grows.  Here one relies on the existence of Euclidean
unimodular lattices $L$ with arbitrarily high $\lambda_1(L)$, as well
as on the (elementary) classification of indefinite odd unimodular
forms, \cf \eqref{81}.  We acknowledge the influence on our approach
of the Lorentzian construction (\ie using indefinite forms) of the
Leech lattice of J.~Conway and N.~Sloane \cite[Chapter 26]{CS}, namely
the following result.

\begin{theorem}
\label{CoSl}
\label{(4.3)} (J. H. Conway, N. J. A. Sloane).  
If
$$t=(3,5,7,\ldots,45, 47, 51 | 145)$$ is a vector with
$q_{24,1}(t)=-1$ in $I_{24,1}$, then $t^\perp\cap I_{24,1}$ is a copy
of the Leech lattice.
\end{theorem}

The {\em first\/} step of our approach can be described as adapting
the Lorentzian construction by replacing the Leech lattice by the
Conway-Thompson lattices.  The latter are unimodular lattices which
define packings of high asymptotic density.  More precisely, we have
the following result \cite[Theorem~9.5]{MH}.

\begin{theorem}[Conway, Thompson]
\label{(1.4)}
For any dimension $n$, there exists a positive definite inner product
space, denoted $\CT_n$, over $\Z$ of odd type and rank $n$ with
\begin{equation}
\label{71}
\min_{x\not=0} x.x\geq k(n),
\end{equation}
where $k(n)$ is asymptotic to $n/2\pi e$ as $n\to \infty$.
\end{theorem}

The {\em second\/} step of our approach is explained in
Section~\ref{(5)}.

\section{Norms
in cohomology and successive minima $\lambda_i$ of lattices}
\label{(two)}

Let $(X,g)$ be a closed orientable Riemannian $(2p)$-dimensional
manifold.  Let $H^p(X,\Z)_\R \subset \H=H^p_{dR}(X)$ be the lattice
defined as the image of $H^p(X,\Z)$ in $H^p(X,\R)$ under the inclusion
$\Z\subset \R$ of coefficients, \ie quotient by its torsion subgroup.
We will sometimes delete the subscript $_\R$, by abuse of notation,
when the torsion subgroup is trivial.  Consider the $L^2$-norm
$|\;|^{\phantom{a}}_{L^2}$ in $\H$, defined by
\begin{equation}
\label{(1.2)}
|f|^2_{L^2}=\int_X f\wedge *f
\end{equation}
for each harmonic $p$-form $f\in \H$, where $*$ is the Hodge operator
for the metric $g$.  The following lemma is obvious, \cf \cite[Lemma
2.21]{FU}.

\begin{lemma}
\label{2.2}
Let $p$ be even.  Then the $L^2$-norm is related to the cup product
form $\omega(f,g)=\int_X f\cup g$ by means of the ``sign reversal''
formula
\begin{equation}
\label{(1.3)}
|f|^2_{L^2} = \langle f,f\rangle= \omega(f^+, f^+) - \omega(f^-, f^-)
\end{equation}
where $f= f^+ + f^-$ is the decomposition
given by the splitting $\H =V^+ + V^-$ into the $(\pm 1)$-eigenspaces
of the involution $*$.
\end{lemma}
Similarly to the notation of formula \eqref{sr} below, we can restate
Lemma~\ref{2.2} as follows:
\begin{equation}
\label{(2.3)}
\langle \;,\;\rangle = \SR(\omega, V^-).
\end{equation}
The lattice $H^p(X,\Z)_\R$ is equipped with the $L^2$-norm defined by
formula~\eqref{(1.2)}.  The dual norm in the similarly defined lattice
$H_p(X,\Z)_\R\subset H_p(X,\R)$ will also be denoted $|\;|_{L^2}$.  

The successive minima are defined as follows.  Note that the second
successive minimum is exploited in Corollary~\ref{(5.1)} below.
\begin{definition}
\label{(2.three)}  
Let $i$ be an integer satisfying $1 \le i\le rk(L)$.  The $i$-th
successive minimum $\lambda_i (L, \|\;\|)$ is the least $\lambda>0$
such that there exist $i$ linearly independent vectors in $L$ of norm
at most $\lambda$:
\[
\lambda_i(L,\|\;\|)= \inf_\lambda \left\{ \lambda\in \R \left|\;
\exists v_1^{\phantom{A}},\ldots,v_i \, (l.i.)\; : \|v_1\|\leq \lambda
,\ldots, \|v_i\|\leq \lambda \right.  \right\} .
\]
\end{definition}

\section{Conformal length and systolic flavors}
\label{(three)}
In this section, we define several flavors of systolic invariants of a
$(2p)$-dimensional Riemannian manifold manifold $(X,g)$.  The (middle
dimensional) conformal $p$-systole, denoted $\confsys_p(g)$, of the
metric $g$, is the least norm of a nonzero element in the integer
lattice in $p$-dimensional cohomology (or, equivalently, homology; see
Remark~\ref{(1.6)}), with respect to the $L^2$-norm \eqref{(1.2)}
defined by $g$:
\[
\begin{aligned}
\confsys_p(g) & =\lambda_1 \left(H^p(X^{2p},\Z)_\R,
|\;|_{L^2}^{\phantom{A}} \right) \\ & = \min\left\{ \left.
|v|^{\phantom{A}}_{L^2} \;\right| v\in{}H^p(X,\Z)_\R \setminus \{0\}
\right\}.
\end{aligned}
\]
Meanwhile, the {\em stable $p$-systole} is the quantity
\[
\stsys_p(g)=\lambda_1(H_p(X^{2p},\Z)_\R,\|\;\|),
\]
where $\|\;\|$ is the stable norm in homology, dual to the comass norm
in cohomology, \cf \cite[4.10]{Fe2}, \cite{BanK}.  The conformal
systole is related to the stable systole as follows:
\begin{equation}
\label{33}
\stsys_p(g)\vol_{2p}(g)^{-\frac{1}{2}}\leq {\binom{2p}{p}}
^{\frac{1}{2}} \confsys_p(g).
\end{equation}
Here the binomial coefficient appears due to the discrepancy between
the linear comass norm and the natural Euclidean norm on the space of
$p$-forms, \cf \cite[section 7]{BanK}.  In the case $p=1$, the
binomial coefficient may be replaced by 1.

\begin{remark}[1-systole asymptotics]
\label{(3.four)}
It should be kept in mind that the asymptotic behavior of the (stable)
1-systole as a function of the genus is completely different from the
conformal systole.  Thus, M.~Gromov \cite[2.C]{Gro2} reveals the
existence of a universal constant $C$ such that we have an
asymptotically {\em vanishing} upper bound
\[
\frac{\sys_1(\Sigma_\gen)^2}{ \area(\Sigma_\gen)} \le C \, \frac{(\log
\gen)^2} {\gen},
\]
for every orientable surface $\Sigma_\gen$ of genus $\gen\ge 2$, with
a Riemannian metric, see \cite[(2.9) and (2.10)]{CK} for related
bounds.  In contrast, P. Buser and P.~Sarnak \cite{BS} provide an
asymptotically {\em infinite\/} lower bound for the maximum of the
conformal systole over the moduli space, \cf inequality \eqref{bs}.
\end{remark}

\begin{remark}[Conformal length]
\label{(1.6)}
The Poincare duality map induces an isometry
\begin{equation}
\label{(3.0)}
\PD:(H^p(X,\Z)_\R,|\;|_{L^2}) \to (H_p(X,\Z)_\R,|\;|_{L^2}),
\end{equation}
proving that the integer lattice in middle dimension is {\em isodual}
in the sense of \cite{CS0,BeM}.  Thus for $p=1$, the invariant
$\confsys_1$ is the conformal length of the surface.
\end{remark}
We have the following upper bound on conformal systole:
\begin{equation}
\label{(3.1)}
\lambda_1(H^p(X^{2p},\Z)_\R, |\;|_{L^2})^2 \leq \gamma_b < {2\over
3}\; b_p(X^{2p}) \;\; \mbox{ for } b_p(X)\geq 2,
\end{equation}
see \cite{BanK} for stable systolic generalisations based on
multiplicative relations in cohomology, and \cite{CK} for an overview.

\section{Systoles of definite intersection forms}
\label{four}
Our main result is Theorem~\ref{T}, which may be viewed as a higher
dimensional generalisation of the Buser-Sarnak theorem \eqref{bs}.  We
briefly discuss the definite case.  Consider the family of manifolds
$n\C P^2$, defined as the connected sum of $n$ copies of the complex
projective plane with the standard orientation.  Recall that these
exhaust the smooth positive definite case by Donaldson's theorem, \cf
\cite{Ka0}.  In contrast to Theorem~\ref{T}, the maximal conformal
systole in the {\em definite\/} case is bounded as the second Betti
number grows:
\begin{equation}
\label{(3.one)}
\lambda_1(H^2(n\C P^2,\Z),|\;|_{L^2}) = \lambda_1(H^2(n\C
P^2,\Z),\sqrt{\omega}) =\lambda_1(\Z^n)= 1,
\end{equation}
for every Riemannian $n\C P^2$, $n=1,2,\ldots$.  This is immediate
from formula \eqref{(1.3)} which identifies the $L^2$-norm and the
intersection form $\omega$ if the latter is positive definite.  By
inequality~\eqref{33}, we obtain the following result, pointed out by
C.~Lebrun: every Riemannian $n\C P^2$ satisfies the inequality
\begin{equation}
\label{42}
\stsys_2\left(n\C P^2 \right)^2\leq 6\;\vol_{4}\left( n\C P^2\right).
\end{equation}

\section{Buser-Sarnak theorem}
\label{(2)}
Our Theorem~\ref{T} may be viewed as a higher dimensional analogue of
the theorem of P. Buser and P. Sarnak \cite[formula (1.13)]{BS}.  Let
$\Sigma_\gen$ be a closed orientable surface of genus $\gen$.  Then
the conformal 1-systole satisfies the bounds
\begin{equation}
\label{bs}
C^{-1}\log \gen < \sup_g \left\{ \lambda_1\left(H^1(\Sigma_\gen,\Z),
|\;|_{L^2}^{\phantom{a}}\right) \right\} ^2 < C \log \gen, \;\forall
\gen\geq 2
\end{equation}
where $C>0$ is a numerical constant, the supremum is over all metrics
$g$ on $\Sigma_\gen$, and $|\;|_{L^2}^{\phantom{a}}$ is the norm
\eqref{(1.2)} associated with $g$.  An explicit upper bound of
$\frac{3}{\pi}\log(4\gen+3)$ is provided in \cite[formula (1.13)]{BS}.

Note that a (weaker) upper bound of $C\sqrt{\gen}$ ( in place of
$C\log \gen$) results from R. Lazarsfeld's work \cite[p.~441,
Proposition, part (i)]{La}.  The systolic quantity $\lambda_1
\left(H^1(\Sigma_\gen,\Z), |\;|_{L^2}^{\phantom{a}}\right)$ may be
viewed as the conformal length of the surface, in view of the
isomorphism of formula~\eqref{(3.0)}.  By conformal invariance, the
supremum in \eqref{bs} may be restricted to the moduli space of
hyperbolic metrics on the surface.

\section{Sign reversal procedure $\SR$ and $Aut(I_{n,1})$-invariance}
\label{six}
Let $q$ be an indefinite quadratic form of index $+1$ (\ie with a
single negative direction) on a vector space $E$ over $\R$, and let
$\v\in E$ be a vector satisfying $q(\v)<0$.  Denote by
$\v^\perp\subset E$ the $q$-orthogonal complement of $\v\in E$, or,
more precisely, the $Q$-orthogonal complement, where
$Q(u,w)=\frac{1}{4} (q(u+w)-q(u-w))$ is the polarisation of $q$.
Thus, we have a decomposition $E=\v^\perp \oplus \R \v $.  The {\em
sign reversal}, $\SR(q,\v)$, is the positive definite form on $E$
obtained by reversing the sign of $q$ in direction $\v$, while keeping
it fixed on $\v^\perp\subset E$:
\begin{equation}
\label{sr}
\SR(q,\v)(x)=q(x^+)-q(x^-),
\end{equation}
where $x=x^+ + x^-$ is the decomposition of $x\in E$ following the
splitting $E= \v^\perp \oplus \R \v $, \cf formula \eqref{(2.3)}.  Let
$\R^{p,q}$ denote the standard real vector space with quadratic form
\begin{equation}
\label{(6.2)}
q_{p,q}(x)=x_1^2+\ldots+x_p^2-x_{p+1}^2-\ldots-x_{p+q}^2,
\end{equation}
and let $I_{p,q}\subset \R^{p,q}$ denote its integer lattice.  For the
purposes of the proof of Theorem~\ref{T}, it is convenient to reverse
the orientation and work instead on the manifold $n\CP^2\#
\overline{\CP}^2$, while hoping that such a step may not prove
baffling to an algebraic geometer.  

Recall that the intersection form on $n\CP^2\# \overline{\CP}^2$ is
$q_{n,1}$, and the integer lattice in two-dimensional homology becomes
a copy of $I_{n,1}$.
\begin{lemma}
\label{(8.1)}
The invariant $\confsys_2\left( n\CP^2\# \overline{\CP}^2,g \right)$
only depends on the orbit of the antiselfdual line of $g$ in
$H^2_{dR}\left( n\CP^2\# \overline{\CP}^2\right)$ under the action of
the automorphism group of $I_{n,1}$.
\end{lemma}

\begin{proof}
An endomorphism $f$ of $H^2_{dR}(n\CP^2\# \overline{\CP}^2)$ which is
an automorphism of the indefinite lattice $I_{n,1}$, induces an
isometry of the definite form $\SR(q_{n,1},\v)$, since $f$ maps the
subspace $\v^{\perp_q}$ to the subspace $f(\v)^{\perp_q}$, and hence
\[
\SR(q_{n,1},\v)(x)=\SR(q_{n,1},f(\v))(f(x)).
\]
Here if $x\in I_{n,1}$, then $f(x)\in I_{n,1}$ by the hypothesis that
$f$ preserves the integer lattice.  Now the lemma follows from the
formula
\[
\confsys_2(g)=\lambda_1\left(H^2 \left( n\CP^2\# \overline{\CP}^2,\Z
\right), \SR\left(q_{n,1}, V^-\right)^{\frac{1}{2}}\right),
\]
where $V^-$ is the antiselfdual direction of $g$.
\end{proof}

\begin{remark}
\label{kotschick}
Note that not all automorphisms of the intersection form can be
realized by a diffeomorphism of the manifold, \cf \cite{Ko}.
\end{remark}

\section{Lorentz construction of Leech lattice and line $\CT_n^\perp$}
\label{8}

Let $\CT_n\subset\R^{n,0}$ be a Conway-Thompson lattice as in
Theorem~\ref{(1.4)}, \ie a unimodular lattice satisfying
$\lambda_1(\CT_n)^2\geq k(n)$.  Then the lattice $\CT_n \oplus
I_{0,1}$ is odd, indefinite, and unimodular, \cf \eqref{(6.2)} and
notation there.  The classification of odd indefinite unimodular forms
\cite[p.~22]{MH} implies that the lattice $I_{n,1}$ contains an
isometric copy of $\CT_n$ such that the $q_{n, 1}$-orthogonal
complement of $\CT_n$ in $I_{n,1}$, is a copy of the line $I_{0,1}$.
In formulas, there exists an isomorphism
\begin{equation}
\label{81}
\phi_n: \CT_n \oplus I_{0,1} \to I_{n,1}
\end{equation}
preserving the bilinear forms.  We will use the following suggestive
notation for the line identified by isomorphism \eqref{81}: let
\begin{equation}
\label{82}
\CT_n^\perp\subset I_{n,1}
\end{equation}
be the $q_{n,1}$-orthogonal complement of $\phi_n(\CT_n\oplus \{0\})
\subset I_{n,1}$.

\begin{corollary}
\label{(5.1)}
Let $I_{n,1}\subset \R^{n,1}$ be the integer lattice.  Let $v\in
I_{n,1}$ be a generator of $\CT_n^\perp\subset I_{n,1}$ as in
\eqref{82}, \ie $v=\phi_n(0,e)$, where $e\in I_{0,1}$ is a generator,
as in isomorphism \eqref{81}.  Consider the norm $
\|x\|_v=\sqrt{\SR(q_{n,1}, v)(x)}, $ in the notation of
formula~\eqref{sr}.  Then the integer lattice has successive minima
$\lambda_1(I_{n,1},\|\;\|_v)=+1$, and
$$
\lambda_2(I_{n,1},\|\;\|_v)^2\geq k(n),
$$
\cf Definition~\ref{(2.three)}, where $k(n)$ is as in
Theorem~\ref{(1.4)}.  In other words, all vectors of square-norm
smaller than $k(n)$ are proportional to each other.
\end{corollary}

\begin{proof}
For any lattice $L$ with a positive definite form, we have the
identity $\SR(L\oplus I_{0,1}, \iota(e)) = L \oplus I_{1,0}$, where
$\iota$ is the inclusion of the second factor.  In particular,
\begin{equation}
\label{(5.2)}
\SR(I_{n,1},\phi_n(\iota(e)))=\CT_n\oplus I_{1,0},
\end{equation}
proving the corollary.
\end{proof}

As an indication of how nontrivial the isomorphism $\phi$ as in
formula~\eqref{81} could be, consider Theorem~\ref{CoSl}, which
exhibits an isomorphism $ \Lambda_{24} \oplus I_{0,1} \to I_{24,1} $,
where $\Lambda_{24}$ is the Leech lattice.

With an eye on the lower bound of our main Theorem~\ref{T}, we first
prove Proposition~\ref{(6.1)} below.  Recall that the intersection
form on $n\CP^2\# \overline{\CP}^2$ is the diagonal form $q_{n,1}$,
\cf~formula~\eqref{(6.2)}.  Let $\phi_n$ be the isomorphism
\eqref{81}.

\begin{proposition}
\label{(6.1)}
If $g$ is a metric on $n\CP^2\# \overline{\CP}^2$ whose antiselfdual
direction is the vector $\phi_n(0,e)\in I_{n,1}$, then all surfaces of
``conformal area'' smaller than $\sqrt{k(n)}$ with respect to $g$ are
homologous to multiples of one another.
\end{proposition}

\begin{proof}
The integer lattice in the selfdual subspace $V^+$ is isometric to the
Conway-Thompson lattice:
\[
V^+ \cap H^2(n\CP^2\# \overline{\CP}^2,\Z)\simeq\CT_n.
\]
Moreover, this copy of the Conway-Thompson lattice is a direct
summand, where the second summand is isometric to $I_{0,1}$.  The sign
reversal formula \eqref{(1.3)} shows that the integer lattice 
\[
\left(H^2( n\CP^2\# \overline{\CP}^2,\Z ),
SR\left( \omega, \CT_n^\perp\right)^{\frac{1}{2}}\right),
\]
is isometric to the positive definite lattice $\CT_n\oplus I_{1,0}$,
where $I_{0,1}$ has been replaced by $I_{1,0}$ as in
formula~\eqref{(5.2)}.  Thus the proposition is a restatement of
Corollary~\ref{(5.1)}.
\end{proof}

\section{Three quadratic forms in the plane}
\label{(5)}

The main result of this section is Lemma~\ref{(4.1)} below on the
interplay of three quadratic forms in the plane, an indefinite one,
$q$, and and a pair of definite ones, $q_1$ and $q_s$, where the
parameter value $s$ will be judiciously chosen in \eqref{(11.12)}.

To go beyond Proposition~\ref{(6.1)} and prove our theorem, the
lattice $\CT_n\oplus I_{1,0}$ is not sufficient, as it contains
vectors of unit norm in the second summand $I_{1,0}$, so that the
quantity $\lambda_1(\CT_n \oplus I_{1,0})=1$ is too small.  In other
words, we need to replace a lower bound for the successive minimum
$\lambda_2$ of the integer lattice, by a lower bound for the
successive minimum $\lambda_1$ for the same lattice, but with respect
to a new norm.  The idea is to deform appropriately the choice of the
negative definite direction $v=\phi_n(0,e)$, responsible for the
Conway-Thompson behavior of its complement.

\begin{figure}
\medskip\noindent
\centerline{{   		%%%the figure
\epsfxsize=6cm
\epsfbox{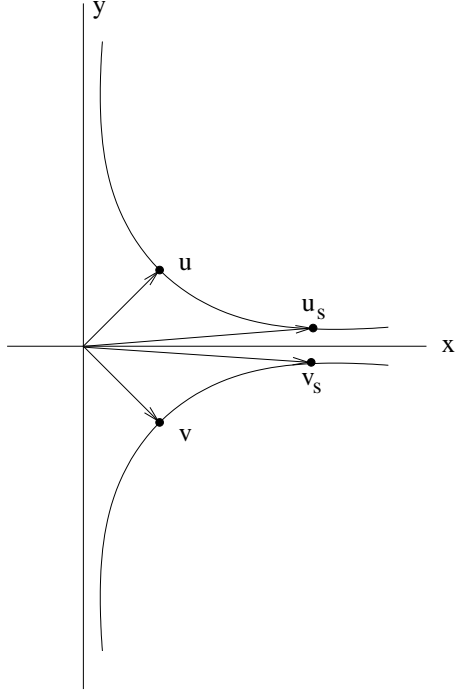}
}}
\caption{Lorentz transformation $A_s$, \cf \eqref{104}}
\label{101}
\end{figure}

Thus, to prove Theorem~\ref{T}, we will apply the surjectivity of the
period map, not to the line $\CT_n^\perp\subset H^2_{dR}\left( n \CP^2
\# \overline{\CP}^2\right) $, but rather to the image of $\CT_n^\perp$
under a suitable ``Lorentz deformation'', \cf Figure~\ref{101} and
formula~\eqref{F}.

\begin{remark}
\label{(3.4)} 
Since the quantity $\lambda_1$ (as in Definition~\ref{(2.three)}) is
continuous as a function on the space of positive definite lattices,
while the form $\SR(\omega ,V)$ is continuous in both parameters, and
$V^-$ depends continuously on the metric, it follows that
Hypothesis~\ref{(1.1)} can be relaxed to assume the density of the
image in place of surjectivity.
\end{remark}

The argument relies on a rather crude bound on the operator norm of
the deformation.  The deformation needs to be sufficient to eliminate
short vectors, but with operator norm controlled so as not to negate
entirely the Conway-Thompson effect.

Sign reversal on the line $\CT_n^\perp\subset I_{n,1}$ produces a
quadratic form with respect to which most vectors are suitably long,
except for a single direction.  To weed out the remaining short
vector, we apply a sutiable deformation, whose linear algebraic
content is presented in Lemma~\ref{(4.1)} below.

Let $\pi$ be the $xy$-plane.  Let $e_1,e_2$ be the standard basis and
$x,y$ the standard coordinates.  Consider the indefinite form
$q=dxdy$, and let $s>0$ be a real parameter.

\begin{definition}
Our ``Lorentz transformation'' $A_s: \pi\to\pi$ is defined by the
matrix $A_s= \begin{pmatrix} s&0\\ 0& \frac{1}{s}\end{pmatrix}$ with
respect to the standard basis, and we set $u_s= A_s (e_1+e_2)=
se_1+{1\over s}e_2$ and $v_s= A_s (e_1-e_2)= se_1-{1\over s}e_2$, as
illustrated in Figure~\ref{101}.
\end{definition}

\begin{lemma}
\label{(4.1)}
Consider the positive definite quadratic form $q_s=\SR(q,v_s)$ on
$\pi$, obtained from $q$ by reversing the sign in the direction $v_s$,
as in formula~\eqref{sr}.  Then the map $A_s:(\pi,q_1)\to (\pi,q_s)$
is an isometry.
\end{lemma}

\begin{proof}  
Since the ``Lorentz transformation'' $A_s$ preserves $q$ and sends
$v_1$ to $v_s$, it is clear that it also sends $q_1$ to $q_s$, but we
will give a short explicit calculation.  We have $q(u_s,v_s)=0$.  Let
$(x',y')$ be the coordinates with respect to the basis $\{u_s,v_s\}$
of the plane $\pi$.  Then the two pairs of coordinates are related by
$x=s(x'+y'),y={1\over s}(x'-y').$ Now $q=dxdy=s(dx'+dy')\,{1\over
s}(dx'-dy')= d{x'}^2-d{y'}^2.$ Therefore by definition,
$q_s=\SR(q,v_s)=d{x'}^2 +d{y'}^2$.  Thus $q_s(u_s,v_s)=0$ and
$q_s(u_s)=q_s(v_s)= 1$.  Similarly, the vectors $u=e_1+e_2$ and
$v=e_1-e_2$ form an orthonormal basis for $q_1$, proving the Lemma.
\end{proof}

\section{Replacing $\lambda_1$ by the geometric mean 
$(\lambda_1\lambda_2)^{1/2}$}
\label{(7)}

Let $L=I_{n,1}\subset (\R^{n,1},q_{n,1})$ be the integer lattice.  Let
$v\in L$ be a vector satisfying $q_{n,1}(v) =-1$ and 
\begin{equation}
\label{(11.1)}
L=\Z v\oplus v^\perp\simeq I_{0,1}\oplus v^\perp,
\end{equation}
where the sublattice 
$\left(v^\perp, \left(q_{n,1} \vert_{v^\perp}
\right)^{\frac{1}{2}} \right)$ is positive definite.  Let
$\SR(q_{n,1},v)$ be the positive definite form obtained by sign
reversal.  Let $\lambda_i= \lambda_i\left(L, \SR(q_{n,1},
v)^{\frac{1}{2}}\right)$ be the successive minima with respect to the
new form.  We have $\lambda_1=1$ but we will ignore this in the
statement of the proposition below, so as to emphasize the geometric
mean inherent in the proof.  Note that
\begin{equation}
\label{(11.2)}
\lambda_2=\lambda_1 \left(v^\perp, \left( q_{n,1} \left| _{v ^\perp}
^{ \phantom {A}} \right.  \right)^{\frac{1}{2}} \right) .
\end{equation}

\begin{proposition}
\label{111}
There is a $q_{n,1}$-preserving transformation $A$ of $\R^{n,1}$ such
that $\lambda_1\left( L, \SR \left( q_{n,1},Av \right) ^{\frac{1}{2}}
\right) \geq \sqrt{\lambda_1\lambda_2}.$
\end{proposition}

\begin{proof}  
Let $\pi\subset\R^{n,1}$ be any 2-plane containing the vector $v$ as
in \eqref{(11.1)}.  We choose coordinates $(x,y)$ in $\pi$ with the
following three properties:
\begin{enumerate}
\item
the union of the $x$-axis and the
$y$-axis in $\pi$ is the intersection of the isotropic cone of
$q_{n,1}$ with $\pi$;
\item
the restriction of $q_{n,1}$ to $\pi$ is the form $q$ of
Lemma~\ref{(4.1)};
\item
with respect to the standard basis $e_1,e_2$ in $\pi$, we have
$v=e_1-e_2$.
\end{enumerate}

Now let $s\in\R$, and set $v_s=se_1-\frac{1}{s} e_2$.  Let $q_s$ be
the positive definite quadratic form obtained by sign reversal
$q_s=\SR(q_{n,1},v_s)$.  Thus, for $s=1$, replacing $q$ by $q_1$ has
the effect of replacing $I_{0,1}$ by $I_{1,0}$ in the
decomposition~\eqref{(11.1)}.  Hence we have the following isometry of
lattices:
\begin{equation}
\label{(7.3)}
(L,q_1)\simeq I_{1,0} \oplus v^\perp.
\end{equation}

We wish to understand the position of the integer lattice $L$ with
respect to the definite form $q_s$ ``deforming'' $q_1$.  By
Lemma~\ref{(4.1)}, the map
\begin{equation}
\label{104}
A_s \oplus {\it Id}_{\pi^\perp},
\end{equation}
also denoted $A_s$, is an isometry from $q_1$ to $q_s$.  Thus the
pullback lattice $\left( A_s^{-1}(L),q_1 \right)$ is isometric to
$(L,q_s)$.  We have $A_s^{-1}(v)=\tfrac{1}{s} e_1 - se_2$, and hence
\begin{equation}
\label{(11.7)}
q_s(v)=q_1\left(A_s^{-1}v\right) = q_1\left(\tfrac{1}{s} e_1 -
se_2\right) = \tfrac{1}{s^2} + s^2 \geq s^2.
\end{equation}
Now consider an element $x\in L=\Z v\oplus v^\perp$ which is not
proportional to the generator $v$ of the first summand.  By the
Pythagorean theorem applied to formula \eqref{(7.3)}, the element $x$
satisfies ${q_1(x)}^{\tfrac{1}{2}}\geq \lambda_1\left(v^\perp\right)=
\lambda_2 \left( L,\sqrt{q_1} \right),$ by formula~\eqref{(11.2)}.
Meanwhile, we have the following bound on the operator norm with
respect to the form $q_1$: $\|A_s\|=\|A_s^{-1}\|\leq s,$ and therefore
\begin{equation}
\label{(11.10)}
q_s(x)= q_1(A_s^{-1} x)\geq {\lambda_2^2\over s^2}.
\end{equation}
Combining \eqref{(11.7)} and \eqref{(11.10)}, we obtain the lower
bound $\lambda_1 \left( L, \sqrt{q_s} \right) \geq \min \left\{ s,
\tfrac {\lambda_2} {s}\right\}$.  Choosing the parameter value
\begin{equation}
\label{(11.12)}
s=\sqrt{\frac{\lambda_2}{\lambda_1}}=\sqrt{\lambda_2},
\end{equation}
we complete the proof
of the proposition.
\end{proof}

\begin{corollary}
\label{(4.three)}  
Let $\Z v =\CT_n^\perp \subset L=\Z^{n,1}$, as in Theorem~\ref{(1.4)}
and isomorphism \eqref{81}.  Then there is a transformation
$A=A_{k(n)^{\frac{1}{4}}}$ of $\R^{n,1}$ such that $\lambda_1\left( L,
\SR ( q_{n,1} ,Av ) ^{\frac{1}{2}}\right)\geq k(n)^{\frac{1}{4}}.$
\end{corollary}

\section{Period map and proof of main theorem}
\label{punchline}

We are now in a position to prove Theorem~\ref{T}.  The
inequality~\eqref{(3.1)} proves the upper bound of
estimate~\eqref{(3.3)}, insofar as $b_2(n\CP \#\overline \CP)=n+1$.
Let us write down a formula, \eqref{F}, for a metric $g_n$ satisfying
the lower bound.  Let $X=n\C P^2\# \overline{\C P^2}$, so that
$H^2(X,\Z)=I_{n,1}$, with cup-form $\omega=q_{n,1}$.  Recall that the
$L^2$-norm in $H^2(X,\R)$ is related to the cup product form
$\omega(f,g)=\int_X f\cup g$ by means of the ``sign reversal'' formula
$ |f|^2_{L^2} = \langle f,f\rangle= \omega(f^+, f^+) - \omega(f^-,
f^-)$, where $f= f^+ + f^-$ is the decomposition given by the
splitting $H^2(X,\R)=V^+ + V^-$ into the $(\pm 1)$-eigenspaces of the
Hodge involution $*$.  It is convenient to introduce the notation \SR,
for the ``sign reversal'' procedure, whose effect is to replace an
indefinite $(n,1)$ form by a positive definite form: $ \langle
\;,\;\rangle = \SR(\omega, V^-), $ \cf formula~\eqref{(2.3)}.

By the Conway-Thompson theorem \cite[Theorem~9.5]{MH}, there exist
positive definite unimodular lattices $\CT_n$ of rank $n$ satisfying
$\lambda_1(\CT_n)^2\geq k(n)$, where $k(n)$ is asymptotic to
$\frac{n}{2\pi e}$ as $n\to \infty$, while $\lambda_1$ is the least
length of a nonzero lattice element, \cf \eqref{lambda}.  Furthermore,
by the classification of odd indefinite unimodular forms
\cite[p.~22]{MH}, there exists a vector $v\in I_{n,1}$ with
$q_{n,1}(v)=-1$, whose orthogonal complement with respect to the
polarisation of $q_{n,1}$ is the lattice $\CT_n$.  Denote by
$\CT_n^\perp\subset H^2(X,\R)$ the negative definite line $\R v$.
Proposition~\ref{111} yields a Lorentzian endomorphism $A_s$ of
$\R^{n,1}$ which replaces the first two successive minima, $\lambda_1$
and $\lambda_2$ (\cf Definition~\ref{(2.three)}), of the lattice with
respect to the definite quadratic form $\SR(\omega,v)$, by their
geometric mean, when one passes to the new definite form
$\SR(\omega,A_s v)$.

Let $\M(X)$ be the space of all Riemannian metrics on $X$, and let
$\G$ be the projectivisation of the negative cone of the form
$\omega$.  Let $\P: \M\to \G$ be the map assigning to each metric, its
antiselfdual direction.  Exploiting the surjectivity of $\P$, we set
\begin{equation}
\label{F}
g_n=\P^{-1} A_{k(n)^{\frac{1}{4}}}\left(\CT_n^\perp\right),
\end{equation}
where $\P^{-1}$ denotes a choice of an inverse image.  Finally, the
lower bound results from the following calculation:
\[
\begin{aligned}
\confsys_2(g_n) &= \lambda_1(H^2(X),|\;|_{L^2}^{\phantom{a}}) \\
&=\lambda_1 \left ( H^2(n\C P^2\# \overline{\C P^2},\Z),
\;\SR\left(q_{n,1}^{\phantom{a}}, A_{k(n)^{\frac{1}{4}}}
\left(\CT_n^\perp\right)\right)^\frac{1}{2} \right) \\ & \geq
\sqrt{\lambda_2 \left( H^2(n\C P^2\# \overline{\C P^2},\Z),
\;\SR\left(q_{n,1}^{\phantom{a}}, \CT_n^\perp\right)^\frac{1}{2}
\right) } \\ & = \sqrt{\lambda_1(\CT_n)} \\ & \geq k(n)^{\frac{1}{4}}.
\end{aligned}
\]

\section*{Acknowledgments}
The author expresses appreciation to P.~Biran, C.~LeBrun, and
S.~Donaldson for insightful comments, and to the referee for a
criticism of an earlier draft that was simultaneously constructive and
exhaustive.

%\vfill\eject
\bibliographystyle{amsalpha}

\end{document}